\documentclass[12pt]{article}
\usepackage[francais,english]{babel}

\usepackage[]{amsfonts}
\usepackage{amssymb}
\usepackage{amsmath}

\def\couleur(#1 #2 #3)
	{
\special{ps::[end]}}

\def\underset#1#2{\mathrel{\mathop{\kern0pt #2}\limits_{#1}}}

\def\overset#1#2{\mathrel{\mathop{\kern0pt #2}\limits^{#1}}}

\def\bx#1{\setbox1=\hbox{\kern3pt{#1}\kern3pt}			
 \dimen1=\ht1 \advance\dimen1 by 3pt \dimen2=\dp1 \advance\dimen2 by 3pt
 \setbox1=\hbox{\vrule height\dimen1 depth\dimen2\box1\vrule}%
 \setbox1=\vbox{\hrule\box1\hrule}%
 \advance\dimen1 by .4pt \ht1=\dimen1
 \advance\dimen2 by .4pt \dp1=\dimen2 \box1\relax}

\def\wbb#1{\kern#1em}
\def\vci{\vrule  width.02em height1.47ex depth-.0ex}		
\def\11{{\rm\wbb{.2}\vci\wbb{-.37}1}}

\newtheorem{Theorem}{Theorem}[section]
\newtheorem{Lemma}[Theorem]{Lemma}
\newtheorem{Remark}[Theorem]{Remark}

\begin{document}

\title{Estimates  $\displaystyle L^{r}-L^{s}$  for solutions of the  $\bar \partial $  equation in strictly pseudo convex domains in  ${\mathbb{C}}^{n}.$ }

\date{}

\author{Eric Amar}
\maketitle
 \ \par 
\ \par 
\ \par 
\ \par 
\ \par 
\renewcommand{\abstractname}{Abstract}

\begin{abstract}
We prove estimates for solutions of the  $\bar \partial u=\omega
 $  equation in a strictly pseudo convex domain  $\displaystyle
 \Omega $  in  ${\mathbb{C}}^{n}.$  For instance if the  $\displaystyle
 (p,q)$  current  $\omega $  has its coefficients  in  $\displaystyle
 L^{r}(\Omega )$  with  $\displaystyle 1\leq r<2(n+1)$  then
 there is a solution  $u$  in  $\displaystyle L^{s}(\Omega )$
  with  $\displaystyle \ \frac{1}{s}=\frac{1}{r}-\frac{1}{2(n+1)}.$
  We also have  $\displaystyle BMO$  and Lipschitz estimates
 for  $\displaystyle r\geq 2(n+1).$  These results were already
 done by S. Krantz~\cite{KrantzDbar76} in the case of  $\displaystyle
 (0,1)$  forms and just for the  $\displaystyle L^{r}-L^{s}$
  part by L. Ma and S. Vassiliadou~\cite{MaVassiliadou00}{\small
   for general  $\displaystyle (p,q)$  forms. To get the complete
 result we propose another approach, based on Carleson measures
 of order  $\alpha $  introduced and studied in~\cite{AmarBonami}
 and on the subordination lemma~\cite{subPrinAmar12}.}\ \par 
\end{abstract}

\section{Introduction.}
\quad \quad  	Let  $\displaystyle \Omega $  be a bounded strictly pseudo
 convex domain with smooth  ${\mathcal{C}}^{\infty }$  boundary.
 We shall denote these domains as s.p.c. domains in the sequel.\ \par 
\quad \quad  	Ovrelid~\cite{Ovrelid71} proved that if we have a  $\displaystyle
 (p,q)$  current  $\omega ,\ \bar \partial $  closed in  $\displaystyle
 \Omega $  and such that its coefficients are in  $\displaystyle
 L^{r}(\Omega )$  then there is a  $\displaystyle (p,q-1)$  current
  $u$  solution of the equation  $\bar \partial u=\omega $  and
 with coefficients still in  $\displaystyle L^{r}(\Omega ).$
  Let us define a norm on these currents :\ \par 
\quad \quad \quad \quad \quad 	 $\displaystyle \omega \in L^{r}_{(p,q)}(\Omega ),\ \omega =\sum_{\left\vert{I}\right\vert
 =p,\left\vert{J}\right\vert =q}{\omega _{I,J}dz^{I}\wedge d\bar
 z^{J}}\Rightarrow {\left\Vert{\omega }\right\Vert}_{r}^{r}:=\sum_{\left\vert{I}\right\vert
 =p,\left\vert{J}\right\vert =q}{{\left\Vert{\omega _{I,J}}\right\Vert}_{r}^{r}}.$
 \ \par 
Then Ovrelid proved that  $\displaystyle \ {\left\Vert{u}\right\Vert}_{r}<C{\left\Vert{\omega
 }\right\Vert}_{r},$  where the constant  $C$  does not depend
 on  $\omega .$ \ \par 
\quad  In the case of  $\displaystyle r=\infty ,$  this was done before
 by Lieb~\cite{Lieb70} and Romanov and Henkin~\cite{RomHenk71}
 proved that still for  $\displaystyle r=\infty ,$  there is
 a solution  $u$  in the space Lipschitz  $\displaystyle 1/2.$
  In the book of Henkin and Leiterer~\cite{HenkinLeiterer84}
 we can find precise references for these topics.\ \par 
\quad \quad  	The  $\displaystyle L^{p}$  results were  strongly improved
 by Krantz~\cite{KrantzDbar76} in the case of  $\displaystyle
 (0,1)$  forms and 	the aim of this work is to generalise Krantz
 results to the case of  $\displaystyle (p,q)$  forms as a consequence
 of results on Carleson measures of order  $\displaystyle \alpha .$ \ \par 
\quad  A more general case was done by  L. Ma and S. Vassiliadou~\cite{MaVassiliadou00}
 on q-convex intersections in  ${\mathbb{C}}^{n},$  but only
 for the  $\displaystyle L^{r}-L^{s}$  part, the Lipschitz one
 is not treated in their work. (Thanks to the referee who signals
 me this nice paper.)\ \par 
\quad \quad  	Moreover, in the case of bounded convex domains of finite type,
 these results are already known, done by K. Diederich, B. Fischer
 and J-E. Fornaess~\cite{DieFisFor99}, A. Cumenge~\cite{Cumeng01}
 and B. Fischer~\cite{Fischer01}.\ \par 
\quad \quad  	So in the case of strictly convex domains, theorem~\ref{LrStrPsConv35}
 can also be seen as a corollary of their results, but for general
 strictly pseudo convex domains this is not the case and of course
 their proofs are much more involved than this one.\ \par 
\quad \quad  	I shall reproof the  $\displaystyle L^{r}-L^{s}$  part of this
  theorem and prove the BMO and Lipschitz one by another approach.\ \par 
\quad \quad  	We already got this kind of results in~\cite{AmarExtForm80}
 by use of Skoda's kernels~\cite{zeroSkoda} but we where dealing
 with boundary values instead of inside ones. Nevertheless using
 Skoda results we shall prove the following theorem, where  $A\lesssim
 B$  means that there is a constant  $\displaystyle C>0$  independent
 of  $A$  and  $B$  such that  $\displaystyle A\leq CB.$ \ \par 
\begin{Theorem}
 ~\label{LrStrPsConv34}Let  $\displaystyle \Omega $  be a s.p.c.
 domain in  ${\mathbb{C}}^{n}$  then for  $\displaystyle 1<r<2n+2$  we have\par 
\quad \quad \quad \quad \quad 	 $\displaystyle \forall \omega \in L^{r}_{(p,q)}(\Omega ),\
 \bar \partial \omega =0,\ \exists u\in L^{s}_{(p,q-1)}(\Omega
 )::\bar \partial u=\omega ,\ {\left\Vert{u}\right\Vert}_{L^{s}(\Omega
 )}\lesssim {\left\Vert{\omega }\right\Vert}_{L^{r}(\Omega )},$ \par 
for any  $s$  such that  $\displaystyle \ \frac{1}{s}>\frac{1}{r}-\frac{1}{2(n+1)}.$
 \par 
\end{Theorem}
\quad \quad  	We shall also generalise Krantz theorem~\cite{KrantzDbar76}
 to  $\displaystyle (p,q)$  forms :\ \par 
\begin{Theorem}
 ~\label{LrStrPsConv35}Let  $\displaystyle \Omega $  be a s.p.c.
 domain in  ${\mathbb{C}}^{n}$  then for  $\displaystyle 1<r<2n+2$  we have\par 
\quad \quad 	 $\bullet $  	 	 $\displaystyle \forall \omega \in L^{r}_{(p,q)}(\Omega
 ),\ \bar \partial \omega =0,\ \exists u\in L^{s}_{(p,q-1)}(\Omega
 )::\bar \partial u=\omega ,\ {\left\Vert{u}\right\Vert}_{L^{s}(\Omega
 )}\lesssim {\left\Vert{\omega }\right\Vert}_{L^{r}(\Omega )},$ \par 
with  $\displaystyle \ \frac{1}{s}=\frac{1}{r}-\frac{1}{2(n+1)}.$ \par 
\quad  $\bullet $ 	 For  $\displaystyle r=2n+2$  we have\par 
\quad \quad \quad \quad \quad   $\displaystyle \exists u\in BMO_{(p,q)}(\Omega )::\bar \partial
 u=\omega ,\ {\left\Vert{u}\right\Vert}_{BMO(\Omega )}\lesssim
 {\left\Vert{\omega }\right\Vert}_{L^{2n+2}(\Omega )}.$ \par 
\quad  If  $\omega $  is a  $\displaystyle (p,1)$  form we have also :\par 
\quad  $\bullet $ 	 for  $\displaystyle r=1,$ \par 
\quad \quad \quad \quad \quad   $\displaystyle \exists u\in L^{s,\infty }_{(p,0)}(\Omega )::\bar
 \partial u=\omega ,\ {\left\Vert{u}\right\Vert}_{L^{s,\infty
 }(\Omega )}\lesssim {\left\Vert{\omega }\right\Vert}_{L^{1}(\Omega )}$ \par 
with  $\displaystyle \ \frac{1}{s}=1-\frac{1}{2(n+1)}.$ \par 
\quad  $\bullet $ 	 for  $\displaystyle r>2n+2,$ \par 
\quad \quad \quad \quad \quad   $\displaystyle \exists u\in \Gamma ^{\beta }_{(p,0)}(\Omega
 )::\bar \partial u=\omega ,\ {\left\Vert{u}\right\Vert}_{\Gamma
 ^{\beta }(\Omega )}\lesssim {\left\Vert{\omega }\right\Vert}_{L^{r}(\Omega
 )},$ \par 
where  $\displaystyle \beta =1-\frac{2n+2}{r}$  and  $\displaystyle
 \ \Gamma ^{\beta }$  is an anisotropic Lipschitz class of functions.\par 
Moreover the solution  $u$  is linear on the data  $\omega .$ \par 
\end{Theorem}
\quad  The classes  $\displaystyle BMO(\Omega )$  and  $\displaystyle
 \Gamma ^{\beta }(\Omega )$  will be defined later. The space
  $\displaystyle L^{s,\infty }_{(p,0)}(\Omega )$  is the Lorentz
 space~\cite{BerghLofstrom76}.\ \par 
\quad  This theorem is 	stronger than theorem~\ref{LrStrPsConv34} because
 here, in the case  $\displaystyle 1\leq r<2(n+1)$  we get the
 result for the end point  $s$  such that  $\displaystyle \ \frac{1}{s}=\frac{1}{r}-\frac{1}{2(n+1)}.$
 \ \par 
\ \par 
Of course if  $\displaystyle u\in L^{s}_{(p,q-1)}(\Omega )$
  for  $\displaystyle s>r$  then  $\displaystyle u\in L^{r}_{(p,q-1)}(\Omega
 )$  hence we also have an strong improvement to Ovrelid's theorem.\ \par 
\quad \quad  	Because the class Lipschitz  $\displaystyle 1/2$  is contained
 in  $\displaystyle \Gamma ^{1}(\Omega )$  we see that we recover
 the Romanov-Henkin result when  $\displaystyle r=\infty $  in
 the case of  $\displaystyle (p,1)$  forms.\ \par 
\ \par 
\quad \quad  	Even if they do not appear in the statement, the Carleson measures
 of order  $\alpha ,$  A. Bonami and I introduced in~\cite{AmarBonami},
 are at the heart of this proof.\ \par 
\ \par 

\section{Proof of the first theorem.}
\quad  Let  $\displaystyle \Omega $  be a s.p.c. in  ${\mathbb{C}}^{n},$
  defined by the function  $\rho \in {\mathcal{C}}^{\infty }({\mathbb{C}}^{n}),$
  i.e.  $\Omega :=\lbrace z\in {\mathbb{C}}^{n}::\rho (z)<0\rbrace
 $  and  $\displaystyle \forall z\in \partial \Omega ,\ \partial
 \rho (z)\neq 0.$ \ \par 
\quad \quad  	Let  $\Omega ':=\lbrace (z,w)\in {\mathbb{C}}^{n}{\times}{\mathbb{C}}::\rho
 '(z,w):=\rho (z)+\left\vert{w}\right\vert ^{2}<0\rbrace $  and
 lift a current  $\omega $  to  $\displaystyle \Omega '$  this
 way :  $\displaystyle \omega '(z,w):=\omega (z).$ \ \par 
\begin{Lemma}
~\label{LpSkoda3} Let  $\displaystyle \Omega $  be a s.p.c.
 domain in  ${\mathbb{C}}^{n},$  with the above notations we have\par 
\quad \quad \quad \quad \quad   $\displaystyle \omega \in L^{r}_{(p,q)}(\Omega )\Rightarrow
 \omega '(z,w)\in L_{(p,q)}^{r}(\partial \Omega ').$ \par 
\end{Lemma}
\quad \quad  	Proof.\ \par 
This is an instance of the subordination principle~\cite{amBerg78},
 ~\cite{subPrinAmar12}. Let  $\displaystyle f(z)\in L^{r}(\Omega
 )$  and set  $\displaystyle f'(z,w):=f(z)$  in  $\displaystyle
 \Omega ',$  then, by the main lemma in~\cite{subPrinAmar12}, p. 6,\ \par 
\quad \quad \quad 	 $\displaystyle \ {\left\Vert{f'}\right\Vert}_{L^{r}(\partial
 \Omega ')}^{r}:=\int_{\partial \Omega '}{\left\vert{f'(z,w)}\right\vert
 ^{r}d\sigma (z,w)}=\int_{\Omega }{\left\vert{f(z)}\right\vert
 ^{r}{\sqrt{-\rho (z)+\frac{\left\vert{{\rm{grad }}\rho (z)}\right\vert
 ^{2}}{4}}}\lbrace \int_{\left\vert{w}\right\vert ^{2}=-\rho
 (z)}{d\left\vert{w}\right\vert }\rbrace dm(z)},$ \ \par 
where  $\displaystyle d\left\vert{w}\right\vert $  is the normalized
 Lebesgue measure~\cite{subPrinAmar12} on the circle  $\displaystyle
 \ \left\vert{w}\right\vert ^{2}=-\rho (z).$  Because  $\displaystyle
 \bar \Omega $  is compact, we have  $\displaystyle \forall z\in
 \bar \Omega \mathrm{,}\mathrm{ }{\sqrt{-\rho (z)+\frac{\left\vert{{\rm{grad
 }}\rho (z)}\right\vert ^{2}}{4}}}\leq C(\rho )<\infty $  hence we have\ \par 
\quad \quad \quad \quad \quad 	 $\displaystyle \ {\left\Vert{f'}\right\Vert}_{L^{r}(\partial
 \Omega ')}^{r}\leq C(\rho )\int_{\Omega }{\left\vert{f(z)}\right\vert
 ^{r}dm(z)}=C(\rho ){\left\Vert{f}\right\Vert}_{L^{r}(\Omega )}.$ \ \par 
It remains to apply this taking for  $f$  any coefficient of
  $\omega .$   $\blacksquare $ \ \par 
\ \par 
\quad \quad  	Proof of theorem~\ref{LrStrPsConv34}.\ \par 
Since  $\displaystyle \Omega $  is a s.p.c. domain so is  $\displaystyle
 \Omega '$  by the subordination lemma~\cite{subPrinAmar12}.
 By use of lemma~\ref{LpSkoda3} we have that  $\displaystyle
 \omega '\in L^{r}_{(p,q)}(\partial \Omega ')$  and still  $\bar
 \partial \omega '=0,$  hence we can apply Skoda's theorem 2
 in~\cite{zeroSkoda} to get that there is a solution  $u'$  of
  $\bar \partial _{b}u'=\omega '$  such that\ \par 
\quad \quad \quad \quad   $\displaystyle u'\in L_{(p,q-1)}^{s}(\partial \Omega ')$  with
  $\displaystyle \ \frac{1}{s}>\frac{1}{r}-\frac{1}{2(n+1)}.$ \ \par 
We have\ \par 
\quad \quad \quad \quad \quad 	 $\displaystyle u'(z,w)=\sum_{I,J}{a'_{I,J}(z,w)dz^{I}\wedge
 d\bar z^{J}}.$ \ \par 
Because  $\omega '$  does not depend on  $w$  we have that the
 coefficients of  $\displaystyle u'$  are holomorphic in  $w,$
  hence we can set (recall that  $\displaystyle u'$  is defined
 on  $\displaystyle \partial \Omega '$ )\ \par 
\quad \quad \quad \quad \quad 	 $\displaystyle \forall z\in \Omega ,\ a_{I,J}(z):=\int_{\left\vert{w}\right\vert
 ^{2}=-\rho (z)}{a'_{I,J}(z,w)d\left\vert{w}\right\vert }$ \ \par 
and\ \par 
\quad \quad \quad \quad \quad 	 $\displaystyle u(z):=\sum_{I,J}{a_{I,J}(z)dz^{I}\wedge d\bar z^{J}},$ \ \par 
then exactly as in~\cite{amExt83} we still have\ \par 
\quad \quad \quad \quad \quad 	 $\bar \partial u=\omega $  in  $\displaystyle \Omega .$ \ \par 
Moreover the subordination lemma~\cite{subPrinAmar12} gives
 again  $\displaystyle u\in L_{(p,q-1)}^{s}(\Omega ),$  because
  $\displaystyle u'\in L_{(p,q-1)}^{s}(\partial \Omega ').$ 
  $\blacksquare $ \ \par 

\section{Carleson measures of order  $\displaystyle \alpha .$ }
\quad \quad  	For  $\displaystyle \Omega $  a s.p.c. domain in  ${\mathbb{C}}^{n},$
  let  $\displaystyle V^{0}(\Omega )$  be the space of bounded
 measures in  $\displaystyle \Omega ,$  and  $\displaystyle V^{1}(\Omega
 )$  the space of Carleson measures in  $\displaystyle \Omega
 $  as defined for instance in~\cite{AmarBonami}. We know that
 these spaces form a interpolating scale for the real method~\cite{AmarBonami},
 and we set\ \par 
\quad \quad \quad \quad \quad 	 $\displaystyle V^{\alpha }(\Omega ):=(V^{0},V^{1})_{(\alpha
 ,\infty )}\ ;\ W^{\alpha }(\Omega ):=(V^{0},V^{1})_{(\alpha
 ,p)}$  with  $\displaystyle p=\frac{1}{1-\alpha }.$ \ \par 
Recall that a  $\displaystyle (p,q)$  form  $\omega $  is in
  $\displaystyle W_{(p,q)}^{\alpha }(\Omega )$  (resp.  $\displaystyle
 V_{(p,q)}^{\alpha }(\Omega )$ ) if its coefficients and the
 coefficients of  $\displaystyle \ \frac{\omega \wedge \bar \partial
 \rho }{{\sqrt{-\rho }}}$  are measures in  $\displaystyle W^{\alpha
 }(\Omega )$  (resp.  $\displaystyle V^{\alpha }(\Omega )$ )
 see~\cite{AmarBonami} and~\cite{AnderCarl00}.\ \par 
A  $\displaystyle (p,q)$  form is in  $\displaystyle L^{r}_{(p,q)}(\Omega
 )$  if just its coefficients are in  $\displaystyle L^{r}(\Omega ).$ \ \par 
\quad \quad  	Let  $\Omega ':=\lbrace (z,w)\in {\mathbb{C}}^{n}{\times}{\mathbb{C}}::\rho
 '(z,w):=\rho (z)+\left\vert{w}\right\vert ^{2}<0\rbrace $  and
 lift a current  $\omega $  to  $\displaystyle \Omega '$  as
 before :  $\displaystyle \omega '(z,w):=\omega (z).$ \ \par 
Our first result links  $\displaystyle L^{r}$  estimates to
 Carleson  $\alpha $  ones.\ \par 
\begin{Theorem}
~\label{LpAmelior42} Let  $\displaystyle \Omega $  be a s.p.c.
 domain in  ${\mathbb{C}}^{n}$ then we have\par 
\quad \quad \quad \quad \quad   $\displaystyle \omega \in L^{r}_{(p,q)}(\Omega )\Rightarrow
 \omega '(z,w):=\omega (z)\in W_{(p,q)}^{\alpha }(\Omega ')$ \par 
with  $\displaystyle \alpha =\frac{1}{r'}+\frac{1}{2(n+1)}.$ \par 
\end{Theorem}
\quad \quad  	Proof.\ \par 
Let  $\displaystyle U':=\bigcup_{j=1}^{N}{Q'(\zeta _{j}',h_{j})\cap
 \partial \Omega '}$  be an open set in  $\displaystyle \partial
 \Omega '$   and  $\displaystyle T(U')=\bigcup_{j=1}^{N}{Q'(\zeta
 _{j}',h_{j})}$  be its associated "tent" set inside~\cite{AmarBonami}
 ; in order to see that a measure  $d\mu =fdm,$  with  $m$  the
 Lebesgue measure in  ${\mathbb{C}}^{n},$  belongs to  $\displaystyle
 V^{\alpha }(\Omega ')$  we have to show, see~\cite{AmarBonami},\ \par 
\quad \quad \quad \quad \quad 	 $\ \int_{T(U')}{\left\vert{f(z')}\right\vert dm(z')}\leq C\left\vert{U'}\right\vert
 ^{\alpha }$ \ \par 
where  $\displaystyle \ \left\vert{U'}\right\vert :=\sigma (U')$
  is the Lebesgue measure of  $\displaystyle U'$  on  $\displaystyle
 \partial \Omega ,$  and with a constant  $C$  independent of
  $\displaystyle U'.$ \ \par 
Because we are dealing with  $\displaystyle (p,q)$  currents
 here, this means that we have to estimate\ \par 
\quad \quad \quad \quad \quad 	 $\displaystyle \ A:=\int_{T(U')}{\frac{\left\vert{\omega (z)}\right\vert
 }{{\sqrt{-\rho '(z,w)}}}dm(z,w)}$ \ \par 
with  $\displaystyle \rho '(z,w):=\rho (z)+\left\vert{w}\right\vert
 ^{2}$  is equivalent to the distance of  $\displaystyle (z,w)\in
 \Omega '$  to the boundary  $\displaystyle \partial \Omega '.$ \ \par 
\quad \quad  	Back to  $A,$ \ \par 
\quad \quad \quad \quad \quad 	 $\displaystyle A:=\int_{T(U')}{\frac{\left\vert{\omega (z)}\right\vert
 }{{\sqrt{-\rho '(z,w)}}}dm(z,w)}\leq \sum_{j=1}^{N}{\int_{Q'_{j}}{\frac{\left\vert{\omega
 (z)}\right\vert }{{\sqrt{-\rho '(z,w)}}}dm(z,w)}}.$ \ \par 
The Carleson window  $\displaystyle Q'_{j}$  is equivalent to
 the product  $\displaystyle (Q'_{j}\cap \partial \Omega '){\times}\lbrack
 h_{j}\rbrack _{\nu _{j}}$  with  $\displaystyle \lbrack h_{j}\rbrack
 _{\nu _{j}}$  the real segment of length  $\displaystyle h_{j}$
  supported by the real normal  $\nu _{j}$  to  $\displaystyle
 \partial \Omega '$  at  $\zeta '_{j}.$  Set  $\displaystyle
 h:=\max _{j=1,...,N}h_{j},$  we shall replace  $\displaystyle
 Q'_{j}$  by  $\displaystyle Q''_{j}:=(Q'_{j}\cap \partial \Omega
 '){\times}\lbrack h\rbrack _{\nu _{j}}.$ \ \par 
\ \par 
So we have\ \par 
\quad \quad \quad \quad \quad 	 $\displaystyle A\leq \sum_{j=1}^{N}{\int_{Q'_{j}}{\frac{\left\vert{\omega
 (z)}\right\vert }{{\sqrt{-\rho '(z,w)}}}dm(z,w)}}\leq \sum_{j=1}^{N}{\int_{Q''_{j}}{\frac{\left\vert{\omega
 (z)}\right\vert }{{\sqrt{-\rho '(z,w)}}}dm(z,w)}},$ \ \par 
where now all the depths have the same value  $\displaystyle
 h.$  Hence by Fubini we have\ \par 
\quad \quad \quad \quad \quad 	 $\displaystyle A\leq \int_{0}^{h}{\frac{1}{{\sqrt{t}}}\lbrace
 \int_{U'_{t}}{\left\vert{\omega (z)}\right\vert d\sigma (z,w)}\rbrace
 dt}$ \ \par 
with   $\displaystyle U'_{t}:=\bigcup_{j=1}^{N}{Q''_{j}\cap
 \partial \Omega '}_{t}$  and  $\displaystyle \partial \Omega
 '_{t}:=\lbrace (z,w)\in \Omega '::\rho (z)+\left\vert{w}\right\vert
 ^{2}=-t\rbrace .$ \ \par 
We can estimate the inner integral by H\"older\ \par 
\quad \quad \quad \quad  	\begin{equation}  \ \int_{U'_{t}}{\left\vert{\omega (z)}\right\vert
 d\sigma (z,w)}\leq {\left({\int_{U'_{t}}{\left\vert{\omega (z)}\right\vert
 ^{r}d\sigma (z,w)}}\right)}^{1/r}{\left({\int_{U'_{t}}{d\sigma
 (z,w)}}\right)}^{1/r'}\label{LrStrPsConv23}\end{equation}\ \par 
but\ \par 
\quad \quad \quad \quad \quad 	 $\displaystyle \ \int_{U'_{t}}{\left\vert{\omega (z)}\right\vert
 ^{r}d\sigma (z,w)}\leq \int_{\partial \Omega _{t}}{\left\vert{\omega
 (z)}\right\vert ^{r}d\sigma (z,w)}\leq C(\rho )\int_{\Omega
 _{t}}{\left\vert{\omega (z)}\right\vert ^{r}\lbrace \int_{\left\vert{w}\right\vert
 ^{2}=-\rho (z)-t}{d\left\vert{w}\right\vert }\rbrace dm(z)}$ \ \par 
where  $\displaystyle d\left\vert{w}\right\vert $  is the normalized
 Lebesgue measure on the circle  $\displaystyle \ \left\vert{w}\right\vert
 ^{2}=-\rho (z)-t.$  Hence, with  $\displaystyle \Omega _{t}:=\lbrace
 z\in \Omega ::\rho (z)<-t\rbrace ,$ \ \par 
\quad \quad \quad \quad \quad 	 $\displaystyle \ \int_{U'_{t}}{\left\vert{\omega (z)}\right\vert
 ^{r}d\sigma (z,w)}\leq C(\rho )\int_{\Omega _{t}}{\left\vert{\omega
 (z)}\right\vert ^{r}dm(z)}=C(\rho ){\left\Vert{\omega }\right\Vert}_{L^{r}(\Omega
 )}^{r}.$ \ \par 
For the last factor of~(\ref{LrStrPsConv23}) we have\ \par 
\quad \quad \quad \quad \quad 	 $\ \int_{U'_{t}}{d\sigma (z,w)}=\sigma (U'_{t})\lesssim \sigma (U'),$ \ \par 
so\ \par 
\quad \quad 	 $\displaystyle A\leq \int_{0}^{h}{\frac{1}{{\sqrt{t}}}\lbrace
 \int_{U'_{t}}{\left\vert{\omega (z)}\right\vert d\sigma (z,w)}\rbrace
 dt}\lesssim {\left\Vert{\omega }\right\Vert}_{L^{r}(\Omega )}(\sigma
 (U'))^{1/r'}\int_{0}^{h}{\frac{dt}{{\sqrt{t}}}}=\frac{1}{2}{\left\Vert{\omega
 }\right\Vert}_{L^{r}(\Omega )}{\sqrt{h}}\sigma (U')^{1/r'}.$ \ \par 
Recall that  $\displaystyle \sigma (Q'_{j})\simeq h_{j}^{(n+1)}$
  then we have\ \par 
\quad \quad \quad \quad \quad 	 $\ {\sqrt{h}}={\sqrt{\max _{j}h_{j}}}\lesssim \max \sigma (Q'_{j})^{1/2(n+1)}\leq
 \sigma (\bigcup_{j=1}^{N}{Q'_{j}\cap \partial \Omega '})^{1/2(n+1)},$ \ \par 
so finally we get\ \par 
\quad \quad \quad \quad \quad 	 $A:=\int_{T(U')}{\frac{\left\vert{\omega (z)}\right\vert }{{\sqrt{-\rho
 '(z,w)}}}dm(z,w)}\lesssim {\left\Vert{\omega }\right\Vert}_{L^{r}(\Omega
 )}\sigma (U')^{\frac{1}{r'}+\frac{1}{2(n+1)}}.$ \ \par 
This means that  $\displaystyle \ \frac{\left\vert{\omega (z)}\right\vert
 }{{\sqrt{-\rho '(z,w)}}}$  is a Carleson measure in  $\displaystyle
 \Omega '$  of order  $\alpha $  with\ \par 
\quad \quad \quad \quad \quad 	 $\displaystyle \alpha =\frac{1}{r'}+\frac{1}{2(n+1)}.$ \ \par 
To get a usual Carleson measure, we need  $\displaystyle \alpha
 =1$  hence\ \par 
\quad \quad \quad \quad \quad 	 $\displaystyle \ \frac{1}{r'}+\frac{1}{2(n+1)}=1\iff r=2(n+1).$ \ \par 
\ \par 
\quad \quad  	We have by theorem 1 in~\cite{AmarBonami}, written in our situation,
 that if  $\displaystyle \mu \in V^{\alpha }(\Omega ')$  then
  $\displaystyle P^{0*}(\mu )\in L^{r,\infty }(\partial \Omega
 '),$  where  $\displaystyle P^{0*}(\mu )$  is the "balayage"
 of  $\mu $  by the Hardy Littlewood kernel  $\displaystyle P_{t}^{0}.$
  Hence we have that the linear operator  $\displaystyle P^{0*}$
  sends  $\displaystyle V^{\alpha _{0}}(\Omega ')$  to  $\displaystyle
 L^{r_{0},\infty }(\partial \Omega '),$  and   $\displaystyle
 V^{\alpha _{1}}(\Omega ')$  to  $\displaystyle L^{r_{1},\infty
 }(\partial \Omega ')$  with, as usual,  $\displaystyle \alpha
 _{j}=1-\frac{1}{r_{j}}.$  This means that\ \par 
\quad \quad \quad \quad \quad 	 $\displaystyle f\in L^{r}(\Omega )\Rightarrow \mu :=f/{\sqrt{-\rho
 '}}dm\in V^{\alpha }(\Omega ')\Rightarrow P^{0*}(\mu )\in L^{s,\infty
 }(\partial \Omega ')$ \ \par 
with control of the norms.\ \par 
So we have a linear operator  $T$  such that, with  $\displaystyle
 r_{0}<r_{1},$ \ \par 
\quad \quad \quad \quad \quad 	 $\displaystyle T\ :\ L^{r_{0}}(\Omega )\rightarrow L^{s_{0},\infty
 }(\partial \Omega '),$  with  $\displaystyle \ \frac{1}{s_{0}}=\frac{1}{r_{0}}-\frac{1}{2(n+1)}\
 ;$ \ \par 
\quad \quad \quad \quad \quad 	 $\displaystyle T\ :\ L^{r_{1}}(\Omega )\rightarrow L^{s_{1},\infty
 }(\partial \Omega '),\ $ with  $\displaystyle \ \frac{1}{s_{1}}=\frac{1}{r_{1}}-\frac{1}{2(n+1)}\
 ;$ \ \par 
hence we can apply Marcinkiewich interpolation theorem between
 these two values of  $r\in \rbrack 1,2(n+1)\lbrack $  i.e.\ \par 
\quad \quad \quad \quad \quad 	 $\displaystyle T\ :\ L^{r}(\Omega )\rightarrow L^{s}(\partial
 \Omega '),$  with  $\displaystyle \ \frac{1}{s}=\frac{1}{r}-\frac{1}{2(n+1)}$
  and  $\displaystyle r\leq s$ \ \par 
which is needed to apply Marcinkiewich theorem, with control of norms.\ \par 
But this implies by theorem 2 in~\cite{AmarBonami}, that  $\displaystyle
 \mu :=f/{\sqrt{-\rho '}}dm\in W^{\alpha }(\Omega ').$   $\blacksquare $ \ \par 

\section{The main result.}
\quad \quad  	Let  $\displaystyle \Omega $  be a domain in  ${\mathbb{C}}^{n}$
  defined by the function  $\rho $  as above ; define  $\Omega
 '\subset {\mathbb{C}}^{n+1}$  the lifted domain : we shall define
 the anisotropic class  $\displaystyle \Gamma ^{\beta }(\partial
 \Omega ')$  as in~\cite{AmarBonami} ; we say that a vector field
  $X$  on  $\displaystyle \partial \Omega '$  is {\sl admissible}
 if  $X$  is of class  ${\mathcal{C}}^{k}$  and at any point
 of  $\displaystyle \zeta \in \partial \Omega ',\ X(\zeta )$
  belongs to the complex tangent space of  $\displaystyle \partial
 \Omega '$  at  $\displaystyle \zeta .$ \ \par 
\quad \quad  	We say that  $\displaystyle u\in \Gamma ^{\beta }(\partial
 \Omega ')$  if  $u$  is bounded on  $\displaystyle \partial
 \Omega '$  and  $u$  belongs to the usual Lipschitz  $\displaystyle
 \Lambda ^{\beta /2}(\partial \Omega '),$  where  $\displaystyle
 \partial \Omega '$  is viewed as a real manifold, and on any
 integral curve of an admissible vector field,  $\displaystyle
 t\in \lbrack 0,1\rbrack \rightarrow \gamma (t)\in \partial \Omega
 ',$  the function  $\displaystyle u\circ \gamma $  belongs to
  $\displaystyle \Lambda ^{\beta }(0,1).$ \ \par 
\quad \quad  	We can now define the class  $\displaystyle \Gamma ^{\beta
 }(\Omega )$  : take a function  $u$  defined in  $\displaystyle
 \Omega $  and lift it as  $\displaystyle u'(z,w):=u(z)$  in
  $\displaystyle \Omega '\ ;$  then  $\displaystyle u\in \Gamma
 ^{\beta }(\Omega )$  if  $\displaystyle u'\in \Gamma ^{\beta
 }(\partial \Omega ').$  	We have that  $\displaystyle u\in \Gamma
 ^{\beta }(\Omega )$  implies that  $\displaystyle u\in L^{\infty
 }(\Omega )$  and  $\displaystyle u\in \Lambda ^{\beta /2}(\Omega
 )$  with a Lipschitz constant uniform in  $\displaystyle \Omega .$ \ \par 
\quad \quad  	The same way we define function  $\displaystyle u\in BMO(\Omega
 )$  if  $\displaystyle u'\in BMO(\partial \Omega ').$  We have
 that  $\displaystyle u\in BMO(\Omega )$  implies that  $\displaystyle
 u\in \bigcap_{r\geq 1}{L^{r}(\Omega )}.$ \ \par 
\quad \quad  	Now we are in position to prove our main result.\ \par 
\begin{Theorem}
 Let  $\displaystyle \Omega $  be a s.p.c. domain in  ${\mathbb{C}}^{n}$
  then for  $\displaystyle 1<r<2n+2$  we have\par 
\quad \quad \quad \quad \quad 	 $\displaystyle \forall \omega \in L^{r}_{(p,q)}(\Omega ),\
 \bar \partial \omega =0,\ \exists u\in L^{s}_{(p,q-1)}(\Omega
 )::\bar \partial u=\omega ,\ {\left\Vert{u}\right\Vert}_{L^{s}(\Omega
 )}\lesssim {\left\Vert{\omega }\right\Vert}_{L^{r}(\Omega )},$ \par 
with  $\displaystyle \ \frac{1}{s}=\frac{1}{r}-\frac{1}{2(n+1)}.$ \par 
For  $\displaystyle r=2n+2$  we have\par 
\quad \quad \quad \quad \quad   $\displaystyle \exists u\in BMO_{(p,q)}(\Omega )::\bar \partial
 u=\omega ,\ {\left\Vert{u}\right\Vert}_{BMO(\Omega )}\lesssim
 {\left\Vert{\omega }\right\Vert}_{L^{2n+2}(\Omega )}.$ \par 
\quad  If  $\omega $  is a  $\displaystyle (p,1)$  form we have also :\par 
for  $\displaystyle r=1,$  we have\par 
\quad \quad \quad \quad \quad   $\displaystyle \exists u\in L^{s,\infty }_{(p,0)}(\Omega )::\bar
 \partial u=\omega ,\ {\left\Vert{u}\right\Vert}_{L^{s,\infty
 }(\Omega )}\lesssim {\left\Vert{\omega }\right\Vert}_{L^{1}(\Omega )}$ \par 
with  $\displaystyle \ \frac{1}{s}=1-\frac{1}{2(n+1)}.$ \par 
for  $\displaystyle r>2n+2$  we have\par 
\quad \quad \quad \quad \quad   $\displaystyle \exists u\in \Gamma ^{\beta }_{(p,0)}(\Omega
 )::\bar \partial u=\omega ,\ {\left\Vert{u}\right\Vert}_{\Gamma
 ^{\beta }(\Omega )}\lesssim {\left\Vert{\omega }\right\Vert}_{L^{r}(\Omega
 )},$ \par 
where  $\displaystyle \beta =1-\frac{2(n+1)}{r}$  and  $\displaystyle
 \ \Gamma ^{\beta }$  is an anisotropic Lipschitz class of functions.\par 
Moreover the solution  $u$  is linear on the data  $\omega .$ \par 
\end{Theorem}
\quad \quad  	Proof.\ \par 
By use of theorem~\ref{LpAmelior42} we have that  $\displaystyle
 \omega '\in W^{\alpha }_{(p,q)}(\Omega ')$  with  $\displaystyle
 \alpha =\frac{1}{r'}+\frac{1}{2(n+1)}$  where  $\displaystyle
 \Omega '$  is still s.p.c.~\cite{subPrinAmar12}, hence we can
 apply the theorem 7 in~\cite{AmarBonami} if  $\omega $  is a
  $\displaystyle (p,1)$  current or the generalisation to  $\displaystyle
 (p,q)$  current done in theorem 4.1 in~\cite{AnderCarl00} to
 get that there is a solution  $u'$  of  $\bar \partial _{b}u'=\omega
 '$  such that\ \par 
\quad \quad \quad \quad   $\displaystyle u'\in L_{(p,q-1)}^{s}(\partial \Omega ')$  with
  $\displaystyle \ \frac{1}{s}=1-\alpha =\frac{1}{r}-\frac{1}{2(n+1)}.$ \ \par 
Because  $\omega '$  does not depend on  $w$  we have that the
 coefficients of  $\displaystyle u'$  are holomorphic in  $w$  hence with\ \par 
\quad \quad \quad \quad \quad 	 $\displaystyle u'(z,w)=\sum_{I,J}{a'_{I,J}(z,w)dz^{I}\wedge
 d\bar z^{J}}$ \ \par 
 we can set (recall that  $\displaystyle u'$  is defined on 
 $\displaystyle \partial \Omega '$ )\ \par 
\quad \quad \quad \quad \quad 	 $\displaystyle \forall z\in \Omega ,\ a_{I,J}(z):=\int_{\left\vert{w}\right\vert
 ^{2}=-\rho (z)}{a'_{I,J}(z,w)d\left\vert{w}\right\vert }$ \ \par 
and we set also\ \par 
\quad \quad \quad \quad \quad 	 $\displaystyle u(z):=\sum_{I,J}{a_{I,J}(z,w)dz^{I}\wedge d\bar
 z^{J}},$ \ \par 
then exactly as in~\cite{amExt83} we still have\ \par 
\quad \quad \quad \quad \quad 	 $\bar \partial u=\omega $  in  $\displaystyle \Omega .$ \ \par 
Moreover the subordination lemma~\cite{subPrinAmar12}, gives
 us  $\displaystyle u\in L_{(p,q-1)}^{s}(\Omega ).$ \ \par 
The last two results came directly from~\cite{AmarBonami}, theorem
 7 and theorem 8 with the fact that we apply them in  $\displaystyle
 \Omega '\subset {\mathbb{C}}^{n+1}$  so we have from theorem
 8 that  $\displaystyle \beta =2(n+1)(\alpha -1).$    $\blacksquare $ \ \par 
\ \par 
\begin{Remark}
 In the range  $\displaystyle 1<r<2n+2$  theorem~\ref{LrStrPsConv35}
 is stronger than theorem~\ref{LrStrPsConv34} because we get
 the result with  $\displaystyle \ \frac{1}{s}=\frac{1}{r}-\frac{1}{2(n+1)}$
  and not only for  $\displaystyle \ \frac{1}{s}>\frac{1}{r}-\frac{1}{2(n+1)}.$
 \par 
\end{Remark}
\ \par 
\ \par 

\bibliographystyle{C:/texlive/2012/texmf-dist/bibtex/bst/base/plain}

\end{document}